# ON STABILITY OF NONZERO SET-POINT FOR NONLINEAR IMPULSIVE CONTROL SYSTEMS

A. D'JORGE, A. L. ANDERSON, A. FERRAMOSCA, A. H. GONZÁLEZ, M. ACTIS

ABSTRACT . The interest in non-linear impulsive systems (NIS) has been growing due to its impact in application problems such as disease treatments (diabetes, HIV, influenza, among many others), where the control action (drug administration) is given by short-duration pulses followed by time periods of null values. Within this framework the concept of equilibrium needs to be extended (redefined) to allows the system to keep orbiting (between two consecutive pulses) in some state space regions out of the origin, according to usual objectives of most real applications. Although such regions can be characterized by means of a discrete-time system obtained by sampling the NIS at the impulsive times, no agreements have reached about their asymptotic stability (AS). This paper studies the asymptotic stability of control equilibrium orbits for NSI, based on the underlying discrete time system, in order to establish the conditions under which the AS for the latter leads to the AS for the former. Furthermore, based on the latter AS characterization, an impulsive Model Predictive Control (i-MPC) that feasibly stabilizes the non-linear impulsive system is presented. Finally, the proposed stable MPC is applied to two control problems of interest: the intravenous bolus administration of Lithium and the administration of antiretrovirals for HIV treatments.

## 1. INTRODUCTION

Hybrid systems (e.g., continuous-time dynamics, discrete-time dynamics, jump phenomena, logic commands, and the like [20]) characterized by abrupt changes in at least one state variable at certain time instants are called Impulsive Systems (IS) [28, 2]. It is well established that an IS provides a natural and reliable framework for mathematical modeling and control of several human diseases [24]. Such results have linked the impulsive behavior with the pills intake (or injection) while the continuous-time counterpart with the pharmacokinetic i.e., absorption and distribution in the organism of such drugs [3]. The dynamics of the Human Immunodeficiency Virus (HIV) has become one of the most studied cases interpreted as an IS since the preliminary result [23]. Other meaningful biomedical problems are Influenza, Ebola [14], Malaria [4], Tumor-bearing [5] and Type I Diabetes [1].

The presence of an impulse effect breaks down fundamental properties of classical dynamical systems necessary to assess the stability [13]. Much of the literature in this regard has focused on the stability of the origin because is the only formal equilibrium state [2, 28, 33], consequently most of the existing control theory assumes the origin as a set-point [15, 6]. However, the control goal in the treatments of the referred biomedical problems is to drive and maintain the system inside a safe zone out of the origin, usually called therapeutic window [31]. Few studies attempted to develop a general theoretical framework for the stability of such a regions. For instance, [30] has linked the stability corresponding to the impulsive time instants with the stability of the IS, however they impose quite restrictive conditions on the control law. Conversely, in [27] a generalized equilibrium set with nonzero set-points was proposed but the result is applicable to the linear case only. To the best of authors knowledge, there is a faltering comprehension on how assessing the stability of regions out of the origin for nonlinear impulsive systems.

---

One of the contribution of this work is to provide a new perspective on stability of nonzero setpoints for constrained IS. The proposed framework allows to prove, by straightforward-methods, that any stabilizer control law of the impulsive time instants can stabilize the free response of the closed-loop impulsive system, and by a set-theoretic analysis it guarantees the feasability of every continuous state between impulses. Besides an useful concept of strong-attractivity, that stands for complementary analysis of several applied problems, is presented. Additionally, we design a nonlinear control predictive strategy - based on the so called zone-MPC (see [10] and references therein) - that feasible stabilize any desired region for the constrained impulsive system. Finally, we illustrate the performance of the strategy by two biomedical applications: the regulation of Lithium ions concentration in a body and the management of HIV treatment in order to maintain an undetectable viral load.

The remaining of the paper is organized as follows: Section 2 introduces some preliminaries concepts about autonomous systems. In Section 3 we present the impulsive control system and the equilibria characterization, Section 4 establishes conditions for the feasibility of the impulsive control system. In Section 5 we present the concepts and results about stability of IS. The control strategy that stabilize the IS is formulated in Section 6. Finally, Section 7 shows some illustrative examples while the concluding remarks are presented in Section 8.

## 2. Preliminaries

### 2.1. Notation.

First we introduce some basic notation. We consider $\mathbb{N}_0$ as the set of non negative integers and $\mathbb{R}^d$ as $d$-dimensional Euclidean space equipped with the euclidean distance between two points defined by $d(x, y) := \|x - y\| = [(x - y)^t(x - y)]^{1/2}$. The euclidean and Hausdorff distance between two sets $X$ and $Y$ of $\mathbb{R}^d$ are given by $d(X, Y) := \inf\{d(x, y) : x \in X, y \in Y\}$ and $d_H(X, Y) = \max\left(\sup_{x \in X} \inf_{y \in Y} d(x, y), \sup_{y \in Y} \inf_{x \in X} d(x, y)\right)$, respectively. The *closed ball with center in* $x \in \mathbb{R}^d$ *and radius* $\varepsilon > 0$ is given by $B(x, \varepsilon) := \{y \in \mathbb{R}^d : d(x, y) \leq \varepsilon\}$. Let $\Omega \subset \mathbb{R}^n$ and $x \in \Omega$, we say that $x$ is an *interior point of* $\Omega$ if the there exist $\varepsilon > 0$ such that $B(x, \varepsilon) \subseteq \Omega$. The *interior of* $\Omega$ is the set of all interior points of $\Omega$ and it is denoted by $\text{int } \Omega$. Finally, given a matrix $A$ of size $m \times n$ we denote (as usual) with $\|A\|_2$ the induced matrix 2-norm, i.e. the square root of the largest eigenvalue of the matrix $A^0 A$ where $A^0$ denotes the conjugate transpose of $A$.

### 2.2. Autonomous systems.

This paper will be based on an autonomous system described by the following set of autonomous nonlinear first-order differential equations

$$\begin{cases} \dot{x}(t) = f(x(t)), & 0 < t < T, \\ x(0) = x_0 \end{cases} \quad (2.1)$$

where the state $x$ takes values in $\mathbb{R}^n$, $t$ represents the continuous-time, $x_0 \in \mathbb{R}^n$ is the initial condition and $T > 0$ is a fixed period of time. As usual to ensure existence and uniqueness of solution of system (2.1), we assume $f$ is Lipschitz continue on a domain $\Omega$. We denote with $\varphi(x_0, t)$ the solution map (sometimes called the *flow*) of system (2.1) at time $t$ with initial condition $x_0$.

**Remark 2.1** (Lipschitz continuity of the solution). *If $f$ is a Lipschitz continuous function on $\Omega$ it can be proved (using the equivalent integral equations and Gronwall's inequality, cf. [7, 32]) that the solution map is locally Lipschitz continuous with respect to the initial condition, i.e. there exists a constant $C$ such that*

$$\|\varphi(x, t) - \varphi(y, t)\| \leq C\|x - y\| \text{ for all } x, y \in \Omega \text{ and } t \in [0, T]. \quad (2.2)$$

*The smallest constant satisfying* (2.2) *is referred as the Lipschitz constant of* $\varphi$ *and is given by*

$$C_\varphi := \sup_{t \in [0,T]} \sup_{x,y \in \Omega} \frac{\|\varphi(x,t) - \varphi(y,t)\|}{\|x - y\|}. \tag{2.3}$$

*This regularity property of the solution will be useful in the proof of the main results of this work.*

## 3. IMPULSIVE CONTROL SYSTEMS

The class of dynamic systems of interest in this paper consists in a set of nonlinear impulsive control system (ICS) of the form

$$\begin{cases} \dot{x}(t) = f(x(t)), & t \neq t_k, \quad k \in \mathbb{N} \\ x(t_k) = x^\circ(t_k) + Bu(t_{k-1}), & k \in \mathbb{N} \\ x(0) = x_0. \end{cases} \tag{3.1}$$

where the independent variable $t \in \mathbb{R}^+$ denotes time, $t_k := kT$, for $k \in \mathbb{N}_0$ and $T > 0$ fixed, denotes the impulse time instants, $T > 0$ the step size in time, $x^\circ(t_k) := \lim_{t \to t_k^-} x(t)$ denotes the state at a time instant $t_k$ right before the impulse, $x(t) \in X \subset \mathbb{R}^n$ denotes the (constrained) state vector, for all $t \geq 0$ and $u(t_k) \in U \subset \mathbb{R}^m$ denotes the (constrained) impulsive controls, for all $k \in \mathbb{N}_0$. Both constraint sets, $X$ and $U$, are assume to be compact convex sets, containing the origin in their interior, while $f(0) = 0$ Matrix $B \in \mathbb{R}^{n \times m}$, which is the impulsive input matrix, is assumed to be full rank.

**Remark 3.1.** *System* (3.1) *is similar to the one reported in* [27] *but with a subtle difference. According to the definition of* $x^\circ(t_k)$, *the discontinuity of the first kind that represents the state jump occurs just before time* $t_k$, *in such a way that* $x(t_k)$ *represents the state after the jump and not before, as it is done in* [27]. *It is a technicality point but it is important since improves the further definitions and simplifies the main results of the paper.*

In case we have a fixed control law $\kappa$ we obtain the following impulsive closed-loop system (ICLS)

$$\begin{cases} \dot{x}(t) = f(x(t)), & t \neq t_k, \quad k \in \mathbb{N} \\ x(t_k) = x^\circ(t_k) + B\kappa(x(t_{k-1})), & k \in \mathbb{N} \\ x(0) = x_0, \end{cases} \tag{3.2}$$

where again more the restrictions are satisfied, i.e. $x(t) \in X$, for all $t \geq 0$, and $\kappa(x(t_k)) \in U$, for all $k \in \mathbb{N}_0$. Note that $x^\circ(t_k) = \varphi(x(t_{k-1}), T)$ hence the above system can be associated to a discrete-time closed-loop system (DCLS)

$$\begin{cases} x(k) = \varphi(x(k-1), T) + B\kappa(x(k-1)), & k \in \mathbb{N} \\ x(0) = x_0, \end{cases} \tag{3.3}$$

where $x(k) \in X_d \subset X$, for all $k \in \mathbb{N}_0$, $\kappa(x(k)) \in U$, for all $k \in \mathbb{N}_0$, and $X_d$ is a compact convex set that contains the origin in its interior. Note that system (3.3) is enough to characterize the impulsive control system (3.2), at sampling times $t_k$. System (3.3) is a sampled-data version of the ICLS (3.2), when the sampling time is given by the impulsive time $T$, i.e. $x(k) := x(t_k)$. This way, the DCLS is the natural way to discretize ICLS, in order to control it by means of conventional discrete-time controllers, such as typical MPCs. In order to satisfy the state constraint of system (3.2) between the jumps, the discrete-time constraint set $X_d$ of system (3.3) must fulfill some properties. Before addressing this point in Section 4, the following fundamental definitions are given.

**Definition 3.2** (Orbit). *Consider the Autonomous System* (2.1). *For any $x \in X$, the orbit of $x$ is given by the set*

$$O_x := \{\varphi(x, \tau) : \tau \in [0, T]\}. \tag{3.4}$$

**Definition 3.3** (Feasibility set). *Consider the ICS* (3.1). *The* feasibility set, $F_X$, *is given by all the points in $X$ whose orbits are also contained in $X$, i.e.*

$$F_X := \{x \in X : O_x \subset X\}. \tag{3.5}$$

**Definition 3.4** (Control equilibria). *Consider the ICS* (3.1). *A state $x_s \in F_X$ is a* control equilibrium point *of the ICS* (3.1) *if there exists an input $u_s = u_s(x_s) \in U$ such that*

$$x_s = \varphi(x_s, T) + Bu_s.$$

*The* control equilibrium set $X_s$ *is the set of every feasible control equilibrium point, i.e.*

$$X_s := \{x \in F_X : \exists u \in U \text{ s.t. } x = \varphi(x, T) + Bu\}.$$

*If $x_s$ is a control equilibrium point for the ICS* (3.1) *we say that $O_s := O_{x_s}$ is a* control equilibrium orbit *for the ICS* (3.1).

The concepts introduced in Definitions 3.2 and 3.3 are generalizations of the concepts presented in [27], in the sense that they refer to any state in the feasible set $X$ and not only to equilibrium states. These generalizations shows to be critical for the developments that follows.

## 4. Feasibility for impulsive control systems

A well known problem of discrete-time systems coming from the sampling of continuous-time ones is how to ensure the feasibility of the solutions of the latter by only imposing constraints on the former [8, 29, 25]. In the case of ICS, this problem is even more difficult to overcome because of the discontinuities in the state trajectories, and it remains still open, although some control strategies were developed to formally tackle the problem [30]. In this work, the feasibility problem (also denoted as the constraints satisfaction problem) is addressed by means of a set-theoretic analysis.

A closed-loop system, DCLS (3.3) (or ICLS (3.2)), is called feasible if for any initial feasible state $x(0) \in X$ all the successive states for DCLS (or all the state trajectories for ICLS), together with the control actions, are also feasible, i.e., $x(t_k) \in X$ for $k \in \mathbb{N}_0$ (or $x(t) \in X$, for $t > 0$), while $\kappa(u(k)) = \kappa(u(t_k)) \in U$ for $k \in \mathbb{N}_0$. So, the following proposition can be stated.

**Proposition 4.1** (Inheritable feasibility). *If $F_X$ is not empty and the constraint set $X_d$ of the DCLS* (3.3) *is such that $X_d \subseteq F_X$ then the feasibility of the DCLS* (3.3) *implies the feasibility of the ICLS* (3.2).

*Proof.* Let $x(t)$ be the solution of the ICLS (3.2). For any $t > 0$ there exist $k \in \mathbb{N}$ and $\tau \in [0, T]$ such that $t = kT + \tau$. Since DCLS (3.3) is feasible then $x(k) \in X_d$ and $u(k) \in U$. By hypothesis $X_d \subseteq F_X$, so $x(k) \in F_X$. From definition (3.5) of $F_X$ we have that $O_{x(k)} \in X$, i.e. $\varphi(x(k), \tau) \in X$, for all $\tau \in [0, T]$. Since $x(t) = \varphi(x(k), \tau)$ for $t = kT + \tau$ the assertion follows.

Note that the non emptiness of $F_X$ depends overall on the size of the impulsive time $T$. If the system is linear and the set polytopic, then there exists an effective method to compute $F_X$ (cf. [18]). However, to the best of our knowledge, for the nonlinear case similar solutions do not exist. Even more in the nonlinear context $F_X$ may not be a convex set. Nevertheless in the following proposition we provide a sufficient condition to guarantee the existence of a constructible compact and convex set $X_d \subseteq F_X$, in the general case.

**Proposition 4.2.** *If there exists $x_* \in X$ such that its orbit $O_*$ (according to Definition (3.2)) satisfies that $O_* \subset \text{int } X$, then there exists a constructible compact convex set $X_d$ such that $X_d \subseteq F_X$.*

*Proof.* For every $t \in [0, T]$ consider $r_t$ as the maximum radius such that a ball centered in $\varphi(x_*, t)$ is contained in $X$, i.e.
$$r_t := \max\{r > 0 : B(\varphi(x_*, t), r) \subseteq X\},$$
and consider $r_*$ as the minimum for all $t \in [0, T]$, i.e.
$$r_* := \min_{t \in [0,T]} r_t.$$
Since $O_* \subset \text{int } X$ it is clear that $r_* > 0$. Let us define
$$X_d := B\left(x_*, \frac{r_*}{C_\varphi}\right) \cap X,$$
where $C_\varphi$ is the Lipschitz constant of $\varphi$ given by 2.1. Note that we need to intersect the ball with $X$ because $C_\varphi$ could be smaller than 1. Since $X_d$ is the intersection of a closed ball with $X$, it is immediate that it is compact and convex. On the other hand, by (2.2) for any $x_d \in X_d$ we have that
$$\|\varphi(x_d, t) - \varphi(x_*, t)\| \leq C_\varphi \|x_d - x_*\| \leq C_\varphi \frac{r_*}{C_\varphi} = r_* \leq r_t,$$
for all $t \in [0, T]$. Then $\varphi(x_d, t) \in B(\varphi(x_*, t), r_t) \subseteq X$, for every $t \in [0, T]$. So $O_{x_d} \subset X$ and by definition (3.5) we get that $x_d \in F_X$. Therefore $X_d \subseteq F_X$, and the assertion follows.

It is worth mentioning that in the examples of Section 7 we will use another strategy for the construction of such a set since, although this approach fully guarantees the inclusion of $X_d$ in $F_X$, it may be somewhat conservative depending on the magnitude of the Lipschitz constant $C_\varphi$. We delve into this in the following remark.

**Remark 4.3.** *In the nonlinear case we can estimate the Lipschitz constant $C_\varphi$ by the following upper bound*
$$C_\varphi \leq e^{T C_f}, \tag{4.1}$$
*where $C_f$ is the Lipschitz constant of $f$ (cf. [32]). If $f$ is differentiable we can bound $C_f$ with the Jacobian of $f$ on $X$, specifically $C_f \leq \max\{\|Jf(x)\|_2 : x \in X\}$. In the linear case (i.e. $f(x)=Ax$) this implies that $C_\varphi \leq e^{T\|A\|_2}$. However, in this case, an upper bound considerably smaller than $e^{T\|A\|_2}$ can be computed, by*
$$C_\varphi \leq \|e^{TA}\|_2 \tag{4.2}$$
*For instance, in the first example presented in Section 7.1 the bound given by (4.1) is $134.37$ while the second one given by (4.2) is $1.13$ (for $T = 3$).*

In any case, it is clear that $C_\varphi \to 1$ when $T \to 0$. Furthermore $r_*$ grows to $r_0$ when $T \to 0$, being $r_0 := \max\{r > 0 : B(x_*, r) \subseteq X\}$. So in the best case $X_d$ is given by the largest ball center in $x_*$ contained in $X$.

## 5. Stability for impulsive systems

The notion of stability for ICLS of nonzero set-points involves some challenges from the theoretical point of view, mostly due to the dual character of this type of systems. Although some preliminary works (cf. [30]) dealt with this issue, they end up being unintuitive or considering quite restrictive assumptions (see for example [30, Proposition 12 and 18]). This is so, mainly, because the stability is defined for arbitrary invariant sets which, in turns, involves intricate definitions in the case of ICS. Our proposal in some sense take a step back it considers first stability of nonzero equilibrium orbits. This is why one of the major contributions of this work is to present a simpler (but nontrivial) approach, meanwhile maintaining its effectiveness both in the theoretical and the practical aspects.

The following definitions are the building blocks to prove the main theorem of the present paper.

### 5.1. New perspective on stability.

From now on let $X_s^*$ be a subset of the control equilibrium set $X_s$ and let us denote by $O_s^*$ the *beam of orbits associated with* $X_s^*$ (see Figure 2), i.e.

$$O_s^* := \bigcup \{O_{x_s} : x_s \in X_s^*\}.$$

**Definition 5.1** (Stability for ICLS). *A beam of orbits $O_s^*$ associated with $X_s^*$ is stable for the impulsive closed-loop system (3.2), if for all $\varepsilon > 0$ there exist $\delta > 0$ such that if $d(x(0), X_s^*) < \delta$ then $d(x(t), O_s^*) < \varepsilon$ for all $t \geq 0$*

**Definition 5.2** (Attractivity for ICLS). *A beam of orbits $O_s^*$ associated with $X_s^*$ is attractive for the impulsive closed-loop system 3.2 if there exists $\Omega \subset X$ such that for every $x_0 \in \Omega$ and every $\varepsilon > 0$ there exists $T = T(x_0, \varepsilon) > 0$ such that $d(x(t), O_s^*) < \varepsilon$ for all $t \geq T$.*

Roughly speaking, what Definition 5.1 states is that: *if a initial state $x_0$ is close enough to a subset $X_s^*$ of the control equilibrium set, then the closed-loop state remains arbitrarily close to the beam of orbits $O_s^*$.* On the other hand, Definition 5.2 states that: *there exists a domain (of attraction) such that if a initial state $x_0$ belongs to such domain then the closed-loop state will eventually be arbitrarily close to the beam of orbits $O_s^*$.*

**Remark 5.3.** *Definitions 5.1 and 5.2 (stability and attractivity for ICLS) are, in some sense, similar to Definition 4 of [30] or Definition 5 of [27] if we consider the sets $Y$ and $Z$ or the sets $X_1$ and $X_2$ as $X_s^*$ and $O_s^*$, respectively. However, our sets are particularized and explicitly linked. This specificity of our definitions translates into better tools to achieve stability proofs. Furthermore we do not require any explicit convexity assumptions on $X_s^*$ and $O_s^*$. Note that even when $X_s^*$ is convex, in general, $O_s^*$ is not.*

**Definition 5.4** (Asymptotic stability for ICLS). *A set of equilibrium orbits $O_s^*$ associated with $X_s^*$ is asymptotically stable for the impulsive closed-loop system (3.2), if it is stable and attractive.*

### 5.2. Stability sufficient conditions.

The following theorem states sufficient conditions for the asymptotic stability of a beam of orbit (even when they do not include the origin), based on the asymptotic stability of a subset of the control equilibrium set associated to the sampled discrete-time system. Definitions of stability and attractivity for discrete-time system can be found in [26].

**Theorem 5.5.** *If a subset of the control equilibrium set $X_s^*$ is asymptotically stable for the discrete-time closed-loop system (3.3), then the beam of orbits $O_s^*$ associated with $X_s^*$ is asymptotically stable for the impulsive closed-loop system (3.2).*

*Proof.* First we will prove stability. Let $\varepsilon > 0$, since $X_s^*$ is stable for the discrete system (3.3), given $\varepsilon_1 = \varepsilon/C_\varphi$ there exists $\delta_1 > 0$ such that

$$\text{if } d(x_0, X_s^*) < \delta_1 \text{ then } d(x(k), X_s^*) < \varepsilon/C_\varphi,$$

for all $k \in \mathbb{N}$. Thus for each $k \in \mathbb{N}$ there exists $x_s^k \in X_s^*$ such that $d(x(k), x_s^k) < \varepsilon/C_\varphi$. Then, considering that $\varphi$ satisfies (2.2), we obtain

$$\|\varphi(x(k), t) - \varphi(x_s^k, t)\| \leq C_\varphi \|x(k) - x_s^k\| < \varepsilon, \text{ for all } t \in [0, T]. \tag{5.1}$$

Hence $d(O_{x(k)}, O_s^*) \leq d(O_{x(k)}, O_{x_s^k}) < \varepsilon$, for all $k \in \mathbb{N}$. Therefore $O_s^*$ is stable for the impulsive closed-loop system (3.2).

Secondly we will prove attractivity. Let $\varepsilon > 0$, since $X_s^*$ is attractive for the discrete system (3.3), then there exists $\Omega \subset X$ such that for every $x_0 \in \Omega$ there exists $K = K(x_0, \varepsilon) \in \mathbb{N}$ such that $d(x(k), X_s^*) < \varepsilon/C_\varphi$, for all $k \geq K$. Following the same argument as before we get that $d(O_{x(k)}, O_s^*) < \varepsilon$, for all $k \geq K$. Therefore $O_s^*$ is attractive for the impulsive closed-loop system (3.2).

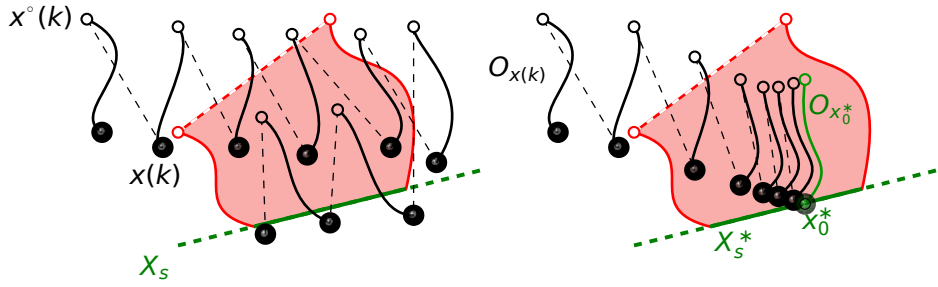

FIGURE 1. Atractivity vs. strong attractivity. On the left we see that the discrete-time closed-loop state $x(k)$ converges to $X_s^*$ but not to any specific point, as a consequence the associated orbits $O_{x(k)}$ do not converge to any specific orbit of the beam of orbits. On the right we observe the opposite, obtaining that $O_{x_0^*}$ is the limit in the sense of the Hausdorff distance of the orbits $O_{x(k)}$.

According to the foregoing result every control strategy capable to steer the DCLS (3.3) to a subset $X_s^*$ of the control equilibrium set, will be able to stabilize the beam of orbits $O_s^*$ associated for the ICLS (3.2). But note that Definition 5.2 (attractivity) requires that the closed-loop state must be close to $O_s^*$, but no qualitative behavior of the trajectory is mentioned, i.e. it could be close but in an erratic way (see Figure 1, left). This situation is sufficient in many cases, such as the biomedical applications mentioned in Section 1 where the goal is just to maintain the system inside a safe zone called therapeutic window (see also the examples treated in Section 7). However a more predictable behavior could be mandatory in some other applications, such as [19], where a *limit orbit* is expected (see Figure 1, right). The following subsection is intended to cover those potentially interesting application cases where a more simple and periodic performance is desired.

5.3. **Further stability definitions for ICLS.** The following definition presents a stronger attractivity concept, that can be shown, as will be stated in Corollary 5.9, when the discrete-time closed-loop system 3.3 converge to an equilibrium point.

**Definition 5.6** (Strong attractivity for ICLS). *A beam of orbits $O_s^*$ associated with $X_s^*$ is strongly attractive for the impulsive closed-loop system* (3.2), *if there exists $\Omega \subset X$ such that for every $x_0 \in \Omega$ there exists an equilibrium point $x_0^* \in X_s^*$ and every $\varepsilon > 0$ there exists $K = K(x_0, \varepsilon) \in \mathbb{N}$ such that $d(\varphi(x(t_k), \tau), \varphi(x_0^*, \tau)) < \varepsilon$, for all $\tau \in [0, T]$ and for all $k \geq K$.*

Note that, as expected, the strong attractivity implies the attractivity of Definition 5.2. Even more, the strong attractivity entails that the orbits $O_{x(t_k)}$ converge uniformly to the (limit) orbit $O_{x_0^*}$ in the sense of the Hausdorff distance (see Figure 1).

**Definition 5.7** (Strong asymptotic stability). *A beam of orbits $O_s^*$ associated with $X_s^*$ is strongly asymptotically stable for the closed-loop system* (3.2), *if it is stable and strongly attractive.*

As stated in the next result, strong asymptotic stability can be easily derived for impulsive closed-loop systems with just the following extra assumption.

**Assumption 5.8.** *For each $x_0$ in the domain of attraction of the discrete-time closed-loop system there exists an equilibrium point $x_0^* \in X_s^*$ to which the closed-loop state converges.*

**Corollary 5.9** (Strong asymptotic stability for ICLS). *If a subset of the control equilibrium set $X_s^*$ is asymptotically stable for the discrete-time closed-loop system* (3.3) *and satisfies Assumption 5.8, then the beam of orbits $O_s^*$ associated with $X_s^*$ is strongly asymptotically stable for the impulsive closed-loop system* (3.2).

*Proof.* The proof follows the same arguments of the proof of Theorem 5.5 replacing $x_s^k$ by $x_0^*$ in (5.1).

**Remark 5.10.** *It is worth mentioning that the definitions and results presented in this section can be easily extended to the case where the subset $X_s^*$ is remaintaind by a general control invariant set and $O_s^*$ by the beam of orbits associated to that invariant. However, since the construction of control invariant sets for impulsive systems is hard and still open problem (see* [30, Section IV-B and C]*) we prefer to focus our work on control equilibrium subsets.*

The next section is focused on the design of a Model Predictive Control that explicitly considers input and state constraints such that the feasibility and stability of the states during the free response is achieved.

## 6. Control Strategy

In the context of stability of nonzero set-points for nonlinear impulsive control systems there are several control strategies that have proven to be useful. Particularly, the predictive control formulation knowing as zone MPC was widely used to stabilize equilibrium regions instead a single set-point (see [12, 11] and the references within it). Moreover, such strategy was previously adjusted for the linear impulsive system framework [27].

In this section a zone MPC formulation is presented to stabilize the nonlinear impulsive system (3.2). Firstly, let us consider a target region not containing the origin, $X^* \subseteq X$, and the following target equilibrium set.

**Definition 6.1** (Target equilibrium set)**.** *Consider the ICLS* (3.2)*. The* target equilibrium set $X_s^*$ *is given by all control equilibrium points inside the target set $X^*$ such that associated bean of orbits belongs to $X^*$, i.e.*

$$X_s^* = \{x_s \in X_s : O_s \subset X^*\} \quad (6.1)$$

Let us consider the following assumption.

**Assumption 6.2.** *The target equilibrium set $X_s^*$ is assumed to be nonempty.*

If Assumption 6.2 is not met, the control problem is not well-posed, meaning that it is not possible - by means of any control strategy- to maintain the impulsive system within the target region $X^*$. To overcome this issue, in [27, Remark 2] a discussion about the relationship between the existence of $X_s^*$ and the length of the selected time period $T$ was made. The study considers a maximal value of $T := T_{max}$, such that Assumption 6.2 holds, and a minimum value of $T := T_{min}$ to represent practical requirements (limitations on the frequency of impulses). In practice we suggest to check if the condition $T_{min} \leq T_{max}$ is fulfilled, otherwise the control problem would not be well-posed.

Consider a control horizon $N \in \mathbb{N}$, a current state $x \in X_d$ and the following cost function

$$J_N(x; \mathbf{u}, x_s, u_s) = \sum_{j=0}^{N-1} \|x_j - x_s\|_Q^2 + \|y - u_s\|_R^2 + \gamma(d_{X_s^*}(x_s) + d_{U_s^*}(u_s)), \quad (6.2)$$

with $\mathbf{u} \in U^N$, $x_s \in X$ and $u_s \in U$ representing the optimization variables. Functions $d_{X_s^*}(x_s)$ and $d_{U_s^*}(u_s)$ stand for the euclidean distance between points and sets, where $U_s^* \subset U$ is the set of admissible inputs necessary to maintain every $x_s \in X_s^*$ in its stationary position. The predicted inputs and states are given by $\mathbf{u} = \{u_0, \cdots, u_{N-1}\}$ and $x_{j+1} = \varphi(x_j, T) + Bu_j$ respectively, with $x_0 = x$ and $j = 0, \cdots, N-1$. The matrices $Q$ and $R$ are positive definite, $\gamma$ is a positive constant, while $x_s$ and $u_s$ are auxiliary reference variables forced to be an equilibrium pair.

Given a time instant $k$ and the current state $x = x(k) \in X_d$, the optimization problem to be solved is given by

$$\min_{\mathbf{u}, x_s, u_s} J_N(x; \mathbf{u}, x_s, u_s) \tag{6.3}$$

$$\text{s.t.} \quad x_0 = x,$$
$$x_{j+1} = \varphi(x_j, T) + Bu(j), j = 0, 1, \cdots, N-1,$$
$$x_j \in X_d, u_j \in U, j = 0, 1, \cdots, N-1,$$
$$x_s \in X, u_s \in U, \tag{6.4}$$
$$x_s = \varphi(x_s, T) + Bu_s,$$
$$x_N = x_s,$$

The control law derived from the application of the receding horizon control policy is given by $\kappa_{\text{MPC}}(x) = u^0(0; x)$ where $u^0(0; x)$ is the first element of the optimal sequence $\mathbf{u}^0(x) = \{u^0(0; x(k)), \ldots, u^0(N-1; x(k))\}$.

In the following lemma, feasibility and stability of the discrete-time closed-loop system are proven.

**Lemma 6.3.** *Consider the discrete-time closed-loop system* (3.3), *with* $\kappa = \kappa_{MPC}$. *Then, the optimization problem* (6.3) *is recursively feasible for all real instant* $k$ *and the control equilibrium set* $X_s^*$ *is asymptotically stable.*

*Proof.* The proof follows the usual steps of asymptotic stability of constrained nonlinear MPC with artificial variables and equality terminal constraint (the details can be seen in the review [17]).

The next result states the asymptotic stability of the impulsive closed-loop nonlinear system (3.2).

**Theorem 6.4.** *Consider the impulsive closed-loop system* (3.2), *with* $\kappa = \kappa_{MPC}$. *The constraint* $x(t) \in X$ *is fulfilled for all* $t \geq 0$ *and the beam of orbits* $O_s^*$ *associated with the target equilibrium set* $X_s^*$ *is asymptotically stable.*

*Proof.* From the constraint $x(k) \in X_d$ for all $k \in \mathbb{N}$ presented on the optimization problem (6.3) along with the result stated on Proposition 4.1, the feasibility of the continuous-time states $x(t) \in X$ follows. On the other hand, Theorem 5.5 and the result on Lemma 6.3 establish the asymptotic stability of the beam of orbits $O_s^*$ associated with $X_s^*$.

**Remark 6.5.** *According to the reports on stability in the zone MPC framework, the closed-loop discrete-time system* (3.3) *converges to the target set* $X_s^*$ *(as was established on Lemma 6.3), but there is no guarantee that it converges to a single point within the target set* $X_s^*$. *So, nothing can be stated about strong asymptotic stability (Definition 5.7) for the impulsive system. However, if we assume that the target set* $X_s^*$ *is a single equilibrium point, the strong stability can be ensured.*

## 7. ILLUSTRATIVE EXAMPLE

In this section, two examples will be presented to test the proposed controller. The first one is a linear system that corresponds to the distribution of lithium in the human body and the second one is a nonlinear system that represents the HIV dynamics.

### 7.1. Lithium Ions Distribution.
In [9] a physiological pharmacokinetic model based on experimental data, which describes the distribution of Lithium ions in the human body upon oral administration, is provided. The system state vector is given by $x(t) = \begin{bmatrix} C_P(t) & C_{RBC}(t) & C_M(t) \end{bmatrix}$, where $C_P(t)$ is the concentration of plasma (P), $C_{RBC}(t)$ is the concentration of the red blood cells, and $C_M(t)$ is the concentration of muscle cells (M). All these concentrations are given in $nmol/L$. The input $u$ is given by the amount of the dose, in $nmol$. The administration period is fixed to $T = 3hr$. The dynamics of the drug distribution is described by

$$\frac{dx(t)}{dt} = \begin{bmatrix} -0.6137 & 0.1835 & 0.2406 \\ 1.2644 & -0.8 & 0 \\ 0.2054 & 0 & -0.19 \end{bmatrix} x(t) \quad (7.1)$$

$$\Delta x(t) = \begin{bmatrix} 10.9 & 0 & 0 \end{bmatrix}^0 u. \quad (7.2)$$

The state and input constraints are given by $X = \{x: \begin{bmatrix} 0 & 0 & 0 \end{bmatrix} \leq x \leq \begin{bmatrix} 2 & 1.2 & 1.2 \end{bmatrix}\}$ and $U = \{u : 0 \leq u \leq 5.95\}$, respectively. The state therapeutic window is defined by $X^* = \{x: \begin{bmatrix} 0.4 & 0.6 & 0.5 \end{bmatrix} \leq x \leq \begin{bmatrix} 0.6 & 0.9 & 0.8 \end{bmatrix}\}$, as it is described in [9] and [30]. The drugs concentration within the boundaries of $X$ guarantees the effectiveness of the therapy. Figure 2 shows the therapeutic window set $X^*$ (green), the control equilibrium target set $X_s^*$ with some target orbits $O_{x_*}$, and discrete feasible set $X_d^*$ (blue). This one is a convex set such that for all $x_0 \in X^*$ it is satisfied that $\varphi(x_0, t) \in X_s^*$ for all $t \in [0, T]$

The MPC controller used is presented in Section 6 and tacking into account [27] is tuned as: $N = 5, Q = diag(1\ 1\ 1), R = 2$, and $\gamma = 100$ In order to observe the dynamic performance of the system, it is considered two starting point $x_0^1 = \begin{bmatrix} 0.2 & 0 & 0 \end{bmatrix}$ and $x_0^2 = \begin{bmatrix} 1.579 & 0 & 0 \end{bmatrix}$.

The results of the simulation are presented in Figures 3 and 5. In particular, Figure 3 shows the state space evolution, the constraint set $X$ (red), the therapeutic window set $X^*$ (green), the control equilibrium target set $X_s^*$, and discrete feasible set $X_d$ (blue).

Concerning the latter, in Proposition 4.1 we present a strategy to construct a compact convex set $X_d$ contained in $F_X$ valid for the general case. As mentioned in Remark 4.3 that approach it may be somewhat conservative depending on the magnitude of the Lipschitz constant $C_\varphi$. Then here we are going to present a simpler construction that remains valid for the linear case due to the convexity of the solution function $\varphi$. Idea of the construction:

1. Mesh the domain $X$ and discretize the time interval $[0, T.]$
2. Select all the points $x$ in the domain mesh that $e^{tA}x$ remain in $X$ for all $t$ in the partition.
3. Make the convex hull of those select points.

The last step it is correct because the solution function $\varphi(x, t) = e^{tA}x$ is convex and $X$ is convex.

It is important to note that the state trajectory converge to $X^*$ as the discrete system tends to a control equilibrium point $x_* \in X_s^*$ determined by the i-NMPC controller. In Figure 4 it is can see a zoom of the graph in the neighborhood of the control equilibrium for the initial state $x_0^1$. Note how by controlling the discrete system, it is directed to a control equilibrium point $x_*$ while the free response converges to the target orbit $O_{x_*}$, represented in red.

Figure 5 shows the state and input time evolution for both initial states. As it is desired, each state is steered to its corresponding therapeutic window. Besides, for the initial state $x_0^1$ the input makes the main effort first and, after its settling time, it remains constant at the desired control equilibrium value. Notice that both, states and inputs, are feasible at any time.

### 7.2. HIV Infection Dynamics With Treatment.

The Human Immunodeficiency Virus (HIV) acts by attacking the immune system, causing its progressive failure over time and its collapse after years (when no treatment is administered).

Several nonlinear models have been developed to describe the dynamics of HIV-1 virus which take into account the kinetics of HIV infection with different cells populations. The virus works by infecting CD4-T cells and its spread can be divided into three stages. The first is called *the acute stage of HIV infection*. In the first days of infection, the virus multiplies rapidly. The spread of the virus activates the immune system to fight the infection. This leads, after a period of 12 to 15 weeks, to the suppression of the spread of the virus and the stabilization of the immune system.

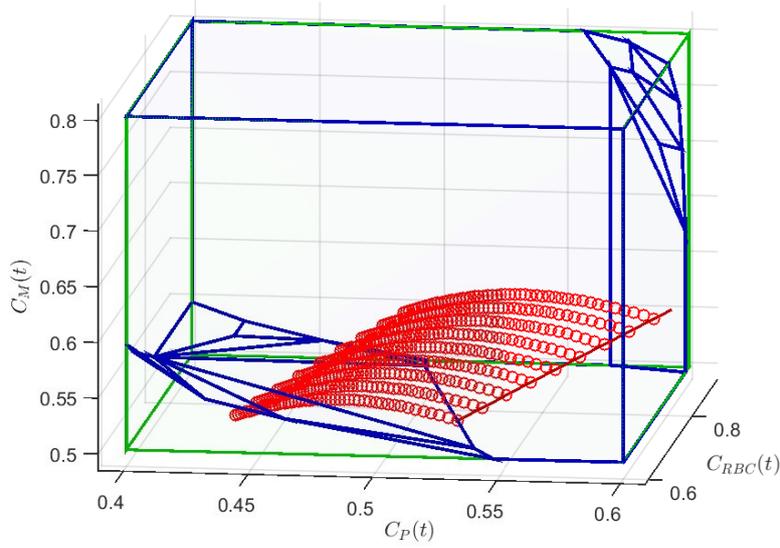

FIGURE 2. Therapeutic window set $X^*$ (green), equilibrium set $X_s^*$ and targets orbit set $O_{x_*}$ (red), and set $X_d^*$ (blue).

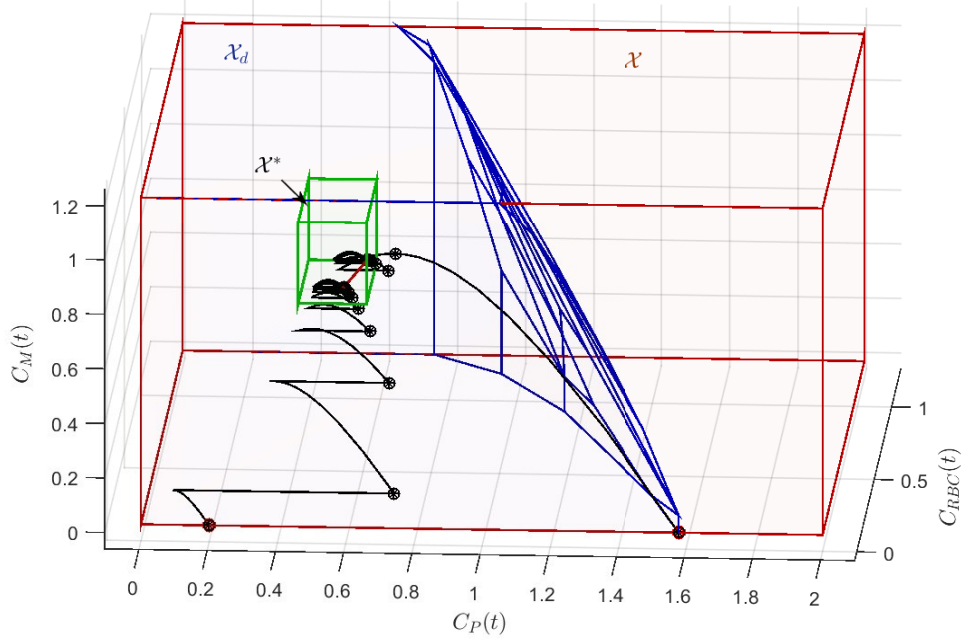

FIGURE 3. Closed-loop system evolution starting from $x_0^1$ and $x_0^2$.

The second is *the clinical latency stage*, also called *chronic HIV infection*. During this stage there is a balance between healthy CD4 + cells and viral load, so the virus is still active but is suppressed by the immune system and reproduces at very low levels. This stage can last up to 10 years for patients not taking medications and up to many decades for patients receiving antiretroviral therapy correctly. Finally, over time, through chronic deterioration, the immune system weakens and becomes vulnerable, which makes the individual vulnerable to opportunistic infections, producing the final stage of HIV infections, the Acquired Immuno Deficiency Syndrome (AIDS).

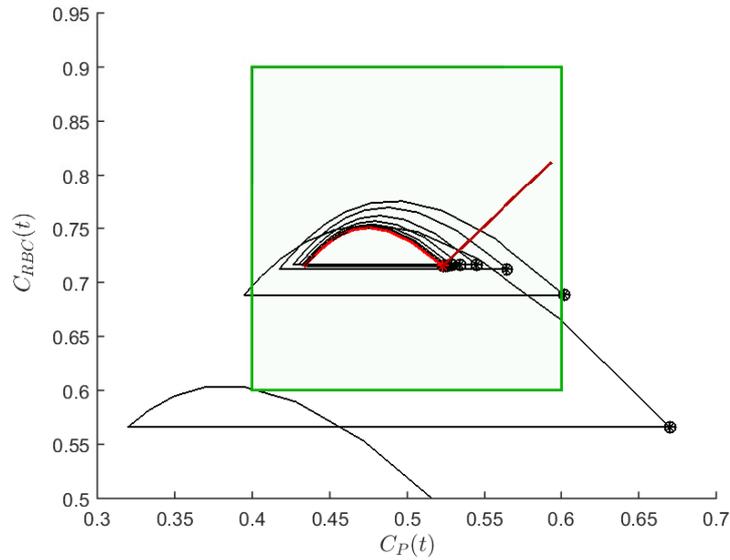

FIGURE 4. Zoom window. It shows the therapeutic window, $X^*$ (green) and the equilibrium set with the target orbit, $O_{x_*}$ (red).

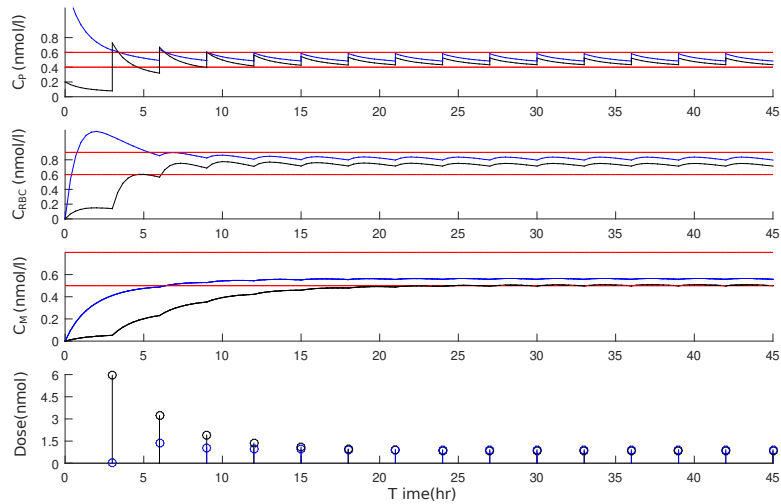

FIGURE 5. Evolution of the states $C_P$, $C_{RBC}$ and $C_M$, and the control input $u$. In black line the time evolution for the initial state $x_0^1$ and in blue line for the initial state $x_0^2$. In red dashed line the limits of the therapeutic window.

Different models can be constructed that describes the effect of HIV on the immune system by considering the interactions between healthy CD4 + T cells, infected CD4 + cells, and viral load. However, for control and parameter estimation based on clinical data purposes, the dynamics of the infection can be modeled by relatively simple ordinary differential equations for the interactions of healthy CD4 + cells ($T_c$), infected CD4 + cells ($y$), free viruses ($z$) [21]. In this paper, the '3-D HIV model' (defined by $T_c$, $y$, $z$) presented in [22] is considered. This model describes the virus infection dynamics and incorporates the interaction of the intake of drugs $(w, u)$ and its concentration in blood according to the notions of pharmacokinetics and pharmacodynamics described in [16]. The complete impulsive model is given by

$$\begin{cases} \dot{T}_c(t) = s - \delta T_c(t) - \beta T_c(t)z(t) \\ \dot{y}(t) = \beta T_c(t)z(t) - \mu y(t) \\ \dot{z}(t) = \left(1 - \frac{w(t)}{w(t) + w_{50}}\right) ky(t) - cz(t) \\ \dot{w}(t) = -K_w w(t) \\ w(\tau_k^+) = w(\tau_k) + u(\tau_k), k \in \mathbb{N} \end{cases}$$

where $T_c$ represents the concentration of healthy CD4 cells [$cell/mm^3$] which are produced from the thymus at a constant rate $s$ [$cell/mm^3 \cdot day$] and die with a half life time equal to $1/\delta$ [$day$]. The healthy cells are infected by the virus at a rate proportional to the product of their population and the amount of free virus particles. Constant $\beta$ [$ml/copies \cdot day$] indicates the effectiveness of the infection process. The infected CD4 cells ($y$) result from the infection of healthy CD4 cells and die at a rate $\mu$ [$1/day$]. Free virus particles ($z$) are produced from infected CD4 cells at a rate $k$ [$copies/cells \cdot mm^3 \cdot ml \cdot day$] and die within a half life time equal to $1/c$ [$day$]. The parameters of the model (taken from [22]) are: $s = 10$, $\delta = 0.02$, $\beta = 2.4 \cdot 10^5$, $\mu = 0.24$, $k = 100$, $c = 2.4$, $K_w = 5.3$ [$day$].

The so-called basic reproduction number (i.e., the coefficient describing the secondary infections produced by an infected cell) is given by $R_0 := \frac{\beta ks}{\mu c \delta}$, and it is assumed to be $R > 1$, to properly describe the spread of the virus in an infected host. Under condition $R > 1$, the control system (7.2) without any treatment (i.e., with $u(t) \equiv 0$) has two equilibrium points, one unstable, given by $x_h := (\frac{s}{\delta}, 0, 0)$ (healthy equilibrium) and one stable, given by $x_e(k) := (\frac{\mu c}{\beta k}, \frac{s}{\mu} - \frac{c\delta}{\beta k}, \frac{sk}{c\mu} - \frac{\delta}{\beta})$ (endemic equilibrium).

The pharmacokinetics and pharmacodynamics phases of the drug administration are related to $w$ (the amount of drug in the body at time $t$) and $\eta = \frac{w(t)}{w(t) + w_{50}}$ (the efficacy of an anti-HIV treatment), where $w_{50}$ is the concentration of drug that lowers the viral load by 50% and in this case we considers $w_{50} = 50$ [$mg$]. Although a cocktail of drugs is generally used, in this simulation only Zidovudine therapies will be considered.

The selected intake period is $T = 0.5$ [$day$]. The state and input constraints are imposed as $X = \{x: [0\ 0\ 0\ 0] \leq x \leq [1200\ 100\ 3000\ 1000]\}$ and $U = \{u: 0 \leq u \leq 610\}$ respectively. The state window target is defined by $X^* = \{x: [900\ 0\ 0\ 0] \leq x \leq [1000\ 5\ 250\ 650]\}$. The control goal is to steer the system from the endemic equilibrium to a healthy zone defined by $X^*$. The anti-HIV treatment is considered successful if both, $z$ is below the threshold of $50$ [$copies/ml$] and $T_c$ is no lower than $900$ [$cell/mm^3$].

The MPC controller is tuned as: $N = 10$, $Q = 5$, $R = 1$, and $\gamma = 5 \cdot 10^6$. In order to test the controller, starting from the initial state $x_0 = [240, 63.33, 2639, 0]$. At this point, the patient is in stage two of the disease, i.e in the clinical latency stage. Treatment begins 20 days after the patient is in this second stage. The results of the simulation are presented in Figures 6 and 7. In particular, Figure 6 shows the state space evolution, the constraint set $X$ (in red), the therapeutic window $X^*$ (green) and the discrete feasible set $X_d$ (blue). As it is desired, the system is steered to therapeutic window. In the right zoom figure, it can be seen in red line part of the equilibrium set and in black solid line how the system tends to the orbit. Figure 7 shows the time evolution of the states $T_c(t)$, $y(t)$, $z(t)$ evolution. On the right, it can be seen a zoom of the graphs which shows how each state is remains into the corresponding therapeutic window.

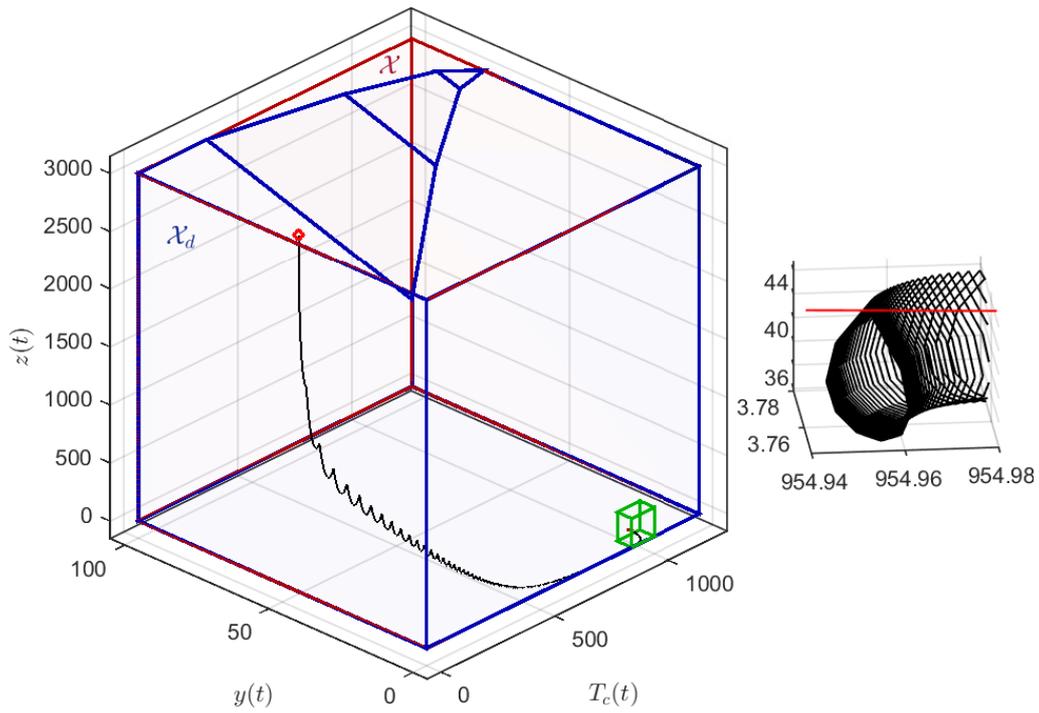

FIGURE 6. Closed-loop system evolution starting from $x_0$. The therapeutic target set is represented by the green box.

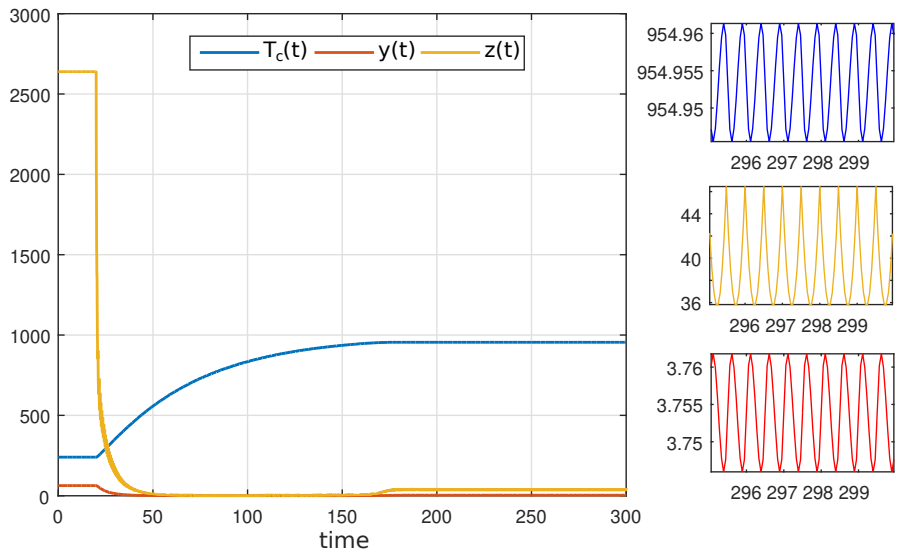

FIGURE 7. Time evolution of the states $T_c(t)$, $y(t)$ and $z(t)$.

## 8. CONCLUSION

In this work we have presented a simpler (but nontrivial) approach to the notion of stability and attractivity of nonzero set-points for impulsive closed-loop system. Furthermore, we also introduced a new concept of strong attractivity that embrace potentially interesting application cases where a smoother performance of the controller is expected near the target. Thanks to all these new definitions we were able to establish sufficient conditions under which the asymptotic stability of a subset of the control equilibrium set associated to the sampled discrete-time system leads to the asymptotic stability for the beam of orbits associated to that subset (even when they do no include the origin).

Based on the latter asymptotic stability characterization, we presented a zone Model Predictive Control (zMPC) that feasibly stabilizes the non-linear impulsive system. The feasibility problem was addressed theoretically providing a minimal condition to guarantee the existence of a constructible set that ensures the feasibility of the solutions of the continuous-time system by only imposing such set as a constraint on the discrete-time system.

Finally, the proposed stable zMPC was applied to two control problems of interest: the intravenous bolus administration of Lithium (linear case) and the administration of antiretrovirals for HIV treatments (nonlinear case). As it was anticipate by the theoretical results, each state was steered to its corresponding therapeutic window.

Future works include an exhaustive study of the construction of control invariant sets for impulsive systems to extend our result to the case where the target set is replaced by this type of sets. They will also include improvements in the algorithms to maintain the feasibility during the free response.